\documentclass[11pt]{article}

\usepackage[T1]{fontenc}
\usepackage[a4paper,margin=1in]{geometry}
\usepackage{amsmath,amssymb,amsthm,mathtools,color}
\usepackage{microtype}
\usepackage{tikz}
\usepackage{enumitem}
\usepackage[hidelinks]{hyperref}
\hypersetup{pdftitle={Sharp minimum-degree thresholds and stability for clique joins of odd cycles},pdfauthor={Jinze Hu, Yanmei Hong, Menghan Ma, Qinghai Liu, Liping Zhang, }}

\allowdisplaybreaks
\newtheorem{theorem}{Theorem}[section]
\newtheorem{lemma}[theorem]{Lemma}

\newtheorem{corollary}[theorem]{Corollary}
\theoremstyle{definition}

\theoremstyle{remark}

\newtheorem*{remark*}{Remark}
\numberwithin{equation}{section}

\newcommand{\N}{\mathbb N}
\newcommand{\chias}{\delta_{\chi}}

\title{Sharp minimum-degree thresholds and stability for clique joins of odd cycles
	\thanks{This work is supported by NSFC (12671401). E-mail: hujinze@fzu.edu.cn (J. Hu), yhong@fzu.edu.cn (Y. Hong), menghanma664@163.com(M. Ma), qliu@fzu.edu.cn (Q. Liu), 1519036064@qq.com (L. Zhang).}
}
\author{Jinze Hu$^a$, Yanmei Hong$^{b,c}$, Menghan Ma$^a$, Qinghai Liu$^{a,c}$\footnote{Corresponding author.}, Liping Zhang$^{a}$ \\[1mm]
	{\small \it
		$^a$Center for Discrete Mathematics, Fuzhou University, Fuzhou, Fujian, 350108 China}\\
	{\small \it
		$^b$School of Mathematics and Statistics, Fuzhou University, Fuzhou, Fujian, 350108 China}\\
	{\small \it
		$^c$Key Laboratory for Operations Research and Cybernetics of Fujian Universities}
}

\date{}
\begin{document}
	\maketitle
	\begin{abstract}
		For an integer $c\ge2$, the chromatic profile of a graph $H$, denoted by $\delta_{\chi}(H,c)$, is the infimum of all $d\in[0,1]$ such that every sufficiently large $H$-free graph $G$ with $\delta(G)\ge d|V(G)|$ is $c$-colourable.
		For vertex-disjoint graphs $H_1$ and $H_2$, let $H_1\vee H_2$ denote their join, obtained by adding all edges between them. 
		We write $C_k$ and $K_r$ for the cycle on $k$ vertices and the complete graph on $r$ vertices, respectively.
		In this paper, we determine $\delta_{\chi}(K_r\vee C_{2k+1},r+2)$ for all integers $r\ge1$ and $k\ge2$. 
		We also establish quantitative stability results describing the structure of non-$(r+2)$-colourable $(K_r\vee C_{2k+1})$-free graphs whose minimum degree is close to the corresponding threshold.
		We further determine the fixed-target homomorphism threshold for odd wheels, proving that $\delta_{\hom}(K_1\vee C_{2k+1},\overline{C_7})=\max\{5/9,(2k+3)/(4k+4)\}$ for every $k\ge2$.
	\end{abstract}

	\section{Introduction}
	All graphs considered in this paper are finite, simple and undirected. A graph is $H$-free if it contains no subgraph isomorphic to $H$; subgraphs are not required to be induced. We write $|G|=|V(G)|$ and $e(G)=|E(G)|$ for the order and size of $G$, and $\delta(G)$ and $\chi(G)$ for its minimum degree and chromatic number. For a positive integer $r$, an $r$-partite graph is a graph with chromatic number at most $r$. Throughout the paper, $\mathbb{N}$ denotes the set of positive integers.
	
	For fixed $r\ge2$, Tur\'an's theorem \cite{Turan} identifies $T_r(n)$, the complete $r$-partite graph with part sizes differing by at most one, as the unique $n$-vertex $K_{r+1}$-free graph with the maximum number of edges. The stability theorem of Simonovits \cite{Sim68} shows that any $n$-vertex $K_{r+1}$-free graph with at least $e(T_r(n))-o(n^2)$ edges can be transformed into $T_r(n)$ by adding or deleting $o(n^2)$ edges. A related problem replaces the edge count by minimum degree and asks for exact colourability. This leads to the chromatic-profile questions considered by Erd\H{o}s and Simonovits \cite{ES73}. For each integer $r\ge2$,  Andr\'asfai, Erd\H{o}s and S\'os \cite{AES} proved that every $K_{r+1}$-free graph of order $n$ with $\delta(G)>\frac{3r-4}{3r-1}n$ is $r$-partite, and the bound is sharp. 
	
	Odd cycles form a basic non-complete family for this problem.
	Andrásfai, Erdős and Sós \cite{AES} proved that every $n$-vertex $\{C_3,C_5,\ldots,C_{2k+1}\}$-free graph with minimum degree greater than $2n/(2k+3)$ is bipartite, and this bound is sharp. When only the single cycle $C_{2k+1}$ is forbidden, Häggkvist \cite{Hagg} established the same conclusion for $k\in\{1,2,3,4\}$, provided that $n$ is sufficiently large.
	For these four lengths, Erd\H{o}s, Faudree, Gy\'arf\'as and Schelp \cite[Theorem~3]{EFGS} identified the balanced blow-up of $C_{2k+3}$ as the only non-bipartite graph attaining equality for sufficiently large $n$. Yuan and Peng \cite{YP24} established the threshold $n/6$ for $C_{2k+1}$-free graphs when $k\ge5$.
	Yuan and Peng \cite{YP} subsequently treated families of odd cycles in a common framework. Let $\mathcal C$ be a family of odd cycles, with longest member $C_{2k+1}$ and shortest omitted odd cycle $C_{2\ell+1}$, where $1\le\ell<k$.  For $n\ge1000k^8$, they proved that every $n$-vertex graph $G$ containing none of the forbidden cycles is bipartite if $\delta(G)>\max\left\{\frac{n}{2(2\ell+1)},\frac{2n}{2k+3}\right\}$,
	and this minimum-degree condition is sharp.
	For $h\in \N$, let $F$ be a graph with vertex set $\{v_1,\ldots,v_h\}$, and let $t\in\mathbb{N}$. The \textit{balanced blow-up }$F(t)$ of $F$ is the graph whose vertex set is partitioned into independent sets $V_1,\ldots,V_h$, each of size $t$, such that, for distinct $i,j$, every vertex of $V_i$ is adjacent to every vertex of $V_j$ if $v_iv_j\in E(F)$, and there are no edges between $V_i$ and $V_j$ otherwise.
	A \textit{blow-up }of $F$ is defined by the same adjacency rule, allowing the independent sets $V_1,\ldots,V_h$ to have different positive sizes.
	For $k\ge 2$, let $\delta_k=\max\{1/6,2/(2k+3)\}$. 
	We present the single-cycle results needed in this paper, using a common order bound.

	\begin{theorem}[\cite{Hagg,YP24,YP}]\label{thm:cycle}
		Let $k\ge2$ and $n\ge1000k^8$ be integers, and let $G$ be an $n$-vertex $C_{2k+1}$-free graph. Then the following hold.
		\begin{enumerate}[label=(\roman*),leftmargin=2em]
			\item  If $\delta(G)>\delta_kn$, then $G$ is bipartite.
			\item  If $k=2$, $G$ is non-bipartite and $\delta(G)\ge2n/7$, then $n=7t$ and $G\cong C_7(t)$ for some $t\in\mathbb N$.
		\end{enumerate}
	\end{theorem}
	
	\begin{remark*}
		Part~(i) follows from H\"aggkvist's result \cite{Hagg} for $2\le k\le4$ and Yuan--Peng's result \cite{YP24} for $k\ge5$. Part~(ii) is the specialization of \cite[Theorem~1.5]{YP} with $\mathcal C=\{C_5\}$ and $\ell=1$.
		For $k\ge5$, the order condition in \cite{YP24} is $n\ge21000k$. Yan, Yuan and Peng \cite[Corollary~1.6]{Yan26} subsequently showed that, for $k\ge79^4$, the weaker condition $n\ge12k+78k^{3/4}$ suffices to guarantee bipartiteness when $\delta(G)>n/6$. In this paper, we use the common order bound stated in Theorem~\ref{thm:cycle}.
	\end{remark*}

	For vertex-disjoint graphs $H_1$ and $H_2$, their join $H_1\vee H_2$ is obtained by adding all edges between them. 
	For $r,k\ge1$, let $H_{r,k}=K_r\vee C_{2k+1}$.
	We write $W_{2k+1}=K_1\vee C_{2k+1}$ for the wheel with $2k+1$ rim vertices and one hub; in particular, $W_3=K_4$. 
	Xu, Jin and Liu \cite{XJL} proved the following tripartiteness result under the additional assumption that $K_4$ is forbidden.
	
	\begin{theorem}[\cite{XJL}]\label{thm:short-wheel}
		Let $G$ be a graph of order $n$ with $\delta(G)>7n/12$. If $G$ is $\{K_4,W_5\}$-free, then $\chi(G)\le3$.
	\end{theorem}

	For $r\ge1$ and $k\ge2$, applying Illingworth's minimum-degree stability theorem \cite[Theorem~15]{IlS} to $H_{r,k}$ gives the sharp threshold for making an $H_{r,k}$-free graph $(r+2)$-partite by deleting $o(n^2)$ edges: it is $(5r+2)/(5r+7)$ when $k=2$ and $(3r+1)/(3r+4)$ when $k\ge3$. More precisely, if $d$ denotes the corresponding value, then, for every fixed $\varepsilon>0$, every $n$-vertex $H_{r,k}$-free graph with minimum degree at least $(d+\varepsilon)n$ can be made $(r+2)$-partite by deleting $o(n^2)$ edges as $n\to\infty$.
	In this paper, we determine the minimum-degree threshold for exact $(r+2)$-colourability, characterize all non-$(r+2)$-colourable $H_{r,k}$-free graphs attaining the threshold, and establish quantitative stability results below it.
	
	Related recent work considers other points of the chromatic profile. B\"ottcher et al.\ \cite{BFMPS} gave three-colourability results for graphs of large odd girth and linear minimum degree. Ning, Wang and Xue \cite{NWX} studied bipartiteness thresholds for three-chromatic forbidden graphs. The forbidden graphs $H_{r,k}$ considered here have chromatic number $r+3\ge4$.

	For a graph $H$ and an integer $c\ge2$, define
	\[
	\begin{aligned}
		\delta_{\chi}(H,c)=\inf\{d\in[0,1]:\ &{\text{there exists }n_0\in\N\text{ such that for every }n\ge n_0,}\\
		&{\text{every }n\text{-vertex }H\text{-free graph }G\text{ with }\delta(G)\ge dn}\\
		&{\text{is }c\text{-colourable}}\}.
	\end{aligned}
	\]
	This parameter fixes the number of colours; it is distinct from the chromatic threshold, which asks only for bounded chromatic number.
	
	By the definition of a balanced blow-up, $K_r(t)$ is the complete $r$-partite graph with each part of size $t$ for $t\in \N$. For convenience, $K_0$ and $K_0(t)$ both denote the graph with no vertices.
	For $r\ge 1$ and $k\ge 2$, let 
	\[
	\delta_{r,k}=\begin{cases}
		\dfrac{5r+2}{5r+7},&k=2,\\[10pt]
		\dfrac{3r+1}{3r+4},&k\ge3.
	\end{cases}                                                  
	\]
	Our first main result is the following.
	\begin{theorem}\label{thm:main}
		For $r\ge1$ and $k\ge2$,
		if $n\ge1000(r+1)k^8$ and $G$ is an $n$-vertex $H_{r,k}$-free graph with $\delta(G)\ge {\color{black}\delta_{r,k}}n$, then either $\chi(G)\le r+2$ \textcolor{black}{or, for some $t\in\mathbb N$, one of the following holds:}
		\begin{enumerate}[label=(\roman*),leftmargin=2em]
			\item $k=2$, $n=(5r+7)t$, and $G\cong K_r(5t)\vee C_7(t)$;
			\item $k\ge3$, $n=(3r+4)t$, and $G\cong K_{r-1}(3t)\vee \overline{C_7}(t)$.
		\end{enumerate}
		Consequently, $\chias(H_{r,k},r+2)={\color{black}\delta_{r,k}}$. 
	\end{theorem}
	
	For sufficiently large graphs, Theorem~\ref{thm:main} strengthens Theorem~\ref{thm:short-wheel} by removing the additional exclusion of $K_4$. Indeed, taking $r=1$ and $k=2$ shows that every sufficiently large $n$-vertex $W_5$-free graph with $\delta(G)>7n/12$ is tripartite. 
	Theorem~\ref{thm:main} also strengthens Illingworth's minimum-degree stability result \cite[Theorem~15]{IlS} for this family by establishing exact $(r+2)$-colourability whenever $\delta(G)>\delta_{r,k}n$. 
	
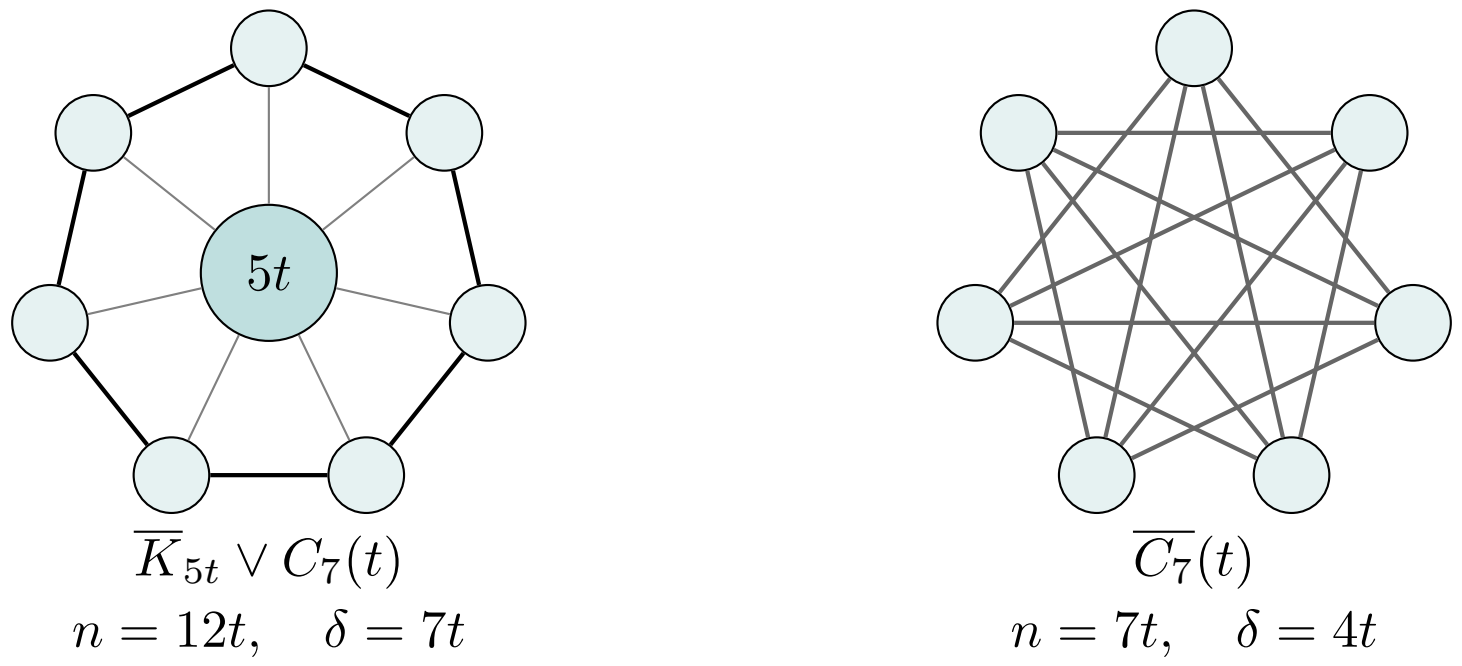
\begin{figure}[htbp]\centering
		\begin{tikzpicture}[scale=.9,every node/.style={font=\small}]
			\begin{scope}[xshift=-3.4cm]
				\foreach \i in {0,...,6}{\node[circle,draw,fill=teal!10,minimum size=5mm,inner sep=0pt](a\i)at({90+360*\i/7}:1.65){};}
				\foreach \i/\j in {0/1,1/2,2/3,3/4,4/5,5/6,6/0}\draw[thick](a\i)--(a\j);
				\node[circle,draw,fill=teal!25,minimum size=9mm](h)at(0,0){$5t$};
				\foreach \i in {0,...,6}\draw[gray](h)--(a\i);
				\node at(0,-2.1){$\overline K_{5t}\vee C_7(t)$};
				\node at(0,-2.65){$n=12t,\quad\delta=7t$};
			\end{scope}
			\begin{scope}[xshift=3.4cm]
				\foreach \i in {0,...,6}{\node[circle,draw,fill=teal!10,minimum size=5mm,inner sep=0pt](b\i)at({90+360*\i/7}:1.65){};}
				\foreach \i/\j in {0/2,0/3,0/4,0/5,1/3,1/4,1/5,1/6,2/4,2/5,2/6,3/5,3/6,4/6}\draw[thick,gray!80!black](b\i)--(b\j);
				\node at(0,-2.1){$\overline{C_7}(t)$};
				\node at(0,-2.65){$n=7t,\quad\delta=4t$};
			\end{scope}
		\end{tikzpicture}
		\caption{\textcolor{black}{The two extremal constructions in Theorem~\ref{thm:main} for $r=1$: $\overline K_{5t}\vee C_7(t)$ on the left and $\overline{C_7}(t)$ on the right. Each small circle represents an independent set of size $t$, and the central circle on the left represents an independent set of size $5t$; lines indicate all edges between the corresponding parts. }}
		\label{fig:examples}
	\end{figure}

	Moreover, Theorem~\ref{thm:main} classifies all non-$(r+2)$-colourable $H_{r,k}$-free graphs of order $n\ge1000(r+1)k^8$ satisfying $\delta(G)=\delta_{r,k}n$, identifying the two extremal families as the only exceptions. 
	\textcolor{black}{Figure~\ref{fig:examples} illustrates the two extremal constructions in Theorem~\ref{thm:main} when $r=1$. For general $r\ge1$, the extremal graphs are obtained by joining the graph on the left with $K_{r-1}(5t)$ and the graph on the right with $K_{r-1}(3t)$, respectively.}
	
	Our second main result describes the non-$(r+2)$-colourable $H_{r,k}$-free graphs of sufficiently large order whose minimum degree is near $\delta_{r,k}n$. We use $E(G)\triangle E(J)$ for the symmetric difference of the two edge sets. 
	
	\begin{theorem}\label{thm:stability-summary}
		Let $r\ge1$, $k\ge2$, $n\ge1000(r+1)k^8$, and
		\[
		0\le\varepsilon<
		\begin{cases}
			1/[(5r+7)(3r+4)],&k=2,\\
			1/[(3r+4)(7r+9)],&k=3,\\
			1/[(3r+4)(4r+5)],&k\ge4.
		\end{cases}
		\]
		Suppose $G$ is an $n$-vertex $H_{r,k}$-free graph with
		$\chi(G)>r+2$ and $\delta(G)\ge({\color{black}\delta_{r,k}}-\varepsilon)n$.
		Then there is a partition of $V(G)$ into nonempty parts with the following properties.
		\begin{enumerate}[label=(\roman*),leftmargin=2em]
			\item If $k=2$, there are $r$ outer parts of size $5n/(5r+7)+O_r(\varepsilon n)$ and seven cycle parts of size $n/(5r+7)+O_r(\varepsilon n)$. If $J$ is the blow-up of $K_r\vee C_7$ on these parts, then
			$|E(G)\triangle E(J)|\le2(5r+7)\varepsilon n^2$.
			\item If $k\ge3$, there are $r-1$ outer parts of size $3n/(3r+4)+O_r(\varepsilon n)$ and seven core parts of size $n/(3r+4)+O_r(\varepsilon n)$. If $J$ is the blow-up of $K_{r-1}\vee \overline{C_7}$ on these parts, then
			$|E(G)\triangle E(J)|
			\le(3r+6)\varepsilon\left(1+\frac{3r+4}{3}\varepsilon\right)n^2$.
			
		\end{enumerate}
	\end{theorem}
	
	For graphs $G$ and $F$, we write $G\to F$ if there exists a graph homomorphism from $G$ to $F$, that is, a map $\varphi:V(G)\to V(F)$ such that $\varphi(u)\varphi(v)\in E(F)$ whenever $uv\in E(G)$.
	For fixed graphs $H$ and $F$, define
	\[
	\begin{aligned}
		\delta_{\hom}(H,F)=\inf\{d\in[0,1]:\ &{\text{there exists }n_0\in\N\text{ such that for every }n\ge n_0,}\\
		&{\text{every }n\text{-vertex }H\text{-free graph }G\text{ with }\delta(G)\ge dn}\\
		&{\text{satisfies }G\to F}\}.
	\end{aligned}
	\]
	Thus $\delta_{\hom}(H,F)$ is the infimum of the minimum-degree bounds that guarantee a homomorphism to the prescribed graph $F$ for every sufficiently large $H$-free graph. The following theorem determines $\delta_{\hom}(W_{2k+1},\overline{C_7})$ for every $k\ge2$.
	
	\begin{theorem}\label{thm:fixed}
		For every integer $k\ge2$, $
		\delta_{\hom}(W_{2k+1},\overline{C_7})=\max\left\{\frac59,\frac{2k+3}{4k+4}\right\}$.
		Moreover, every $n$-vertex $W_{2k+1}$-free graph $G$ with $n\ge2000k^8$ and $\delta(G)>\delta_{\hom}(W_{2k+1},\overline{C_7})n$ admits a homomorphism to $\overline{C_7}$.
	\end{theorem}
	For $k\ge3$, the threshold for a homomorphism to $\overline{C_7}$ is strictly below the tripartiteness threshold $4/7$ for $W_{2k+1}$-free graphs.
	
	In Section~\ref{sec:thresholds}, we prove Theorem~\ref{thm:main}, determining the sharp colourability thresholds for $H_{r,k}$-free graphs and classifying the extremal graphs. Section~\ref{sec:fixed-target} gives the proof of Theorem~\ref{thm:fixed}, which determines the minimum-degree threshold for a $W_{2k+1}$-free graph to admit a homomorphism to $\overline{C_7}$. Section~\ref{sec:stability} proves the stability result stated in Theorem~\ref{thm:stability-summary}. Theorem~\ref{thm:short-stability} gives the case $k=2$, and we deduce the case $k\ge3$ from Theorem~\ref{thm:local-stability} on $r$-locally bipartite graphs.

	\section{The sharp colourability threshold}\label{sec:thresholds}
	
	For $R\subseteq V(G)$, let $N_G(R)=\bigcap_{v\in R}N_G(v)$. We identify a clique with its vertex set when writing its common neighbourhood. A graph $G$ is \emph{$r$-locally bipartite} if $G[N_G(R)]$ is bipartite for every $r$-clique $R$. For $r=1$, this is ordinary \textit{local bipartiteness}.
	The proof of Theorem~\ref{thm:main} treats $k\ge3$ and $k=2$ separately. For $k\ge3$, the minimum-degree assumption and Theorem~\ref{thm:cycle} imply that the common neighbourhood of every $r$-clique induces a bipartite graph. We can therefore apply Illingworth's results \cite{IlB,IlC} to determine the colourability and equality cases. For $k=2$, if $G$ is not $(r+2)$-colourable, we find an $r$-clique whose common neighbourhood induces a non-bipartite graph. Theorem~\ref{thm:cycle}(ii) identifies this graph as $C_7(t)$, and the minimum-degree condition then forces the remaining vertices to induce $K_r(5t)$, with all edges present between these two graphs.
	We use Illingworth's local colourability theorem \cite[Theorem~5.6(a),(b)]{IlC} and his homomorphism theorem and sharpness examples for locally bipartite graphs \cite[Theorem~2 and Fig.~2]{IlB}.
	\begin{theorem}[\cite{IlB,IlC}]\label{thm:local}
		Let $r\ge1$ be an integer and let $G$ be an $r$-locally bipartite graph.
		\begin{enumerate}[label=(\alph*),leftmargin=2em]
			\item  If $\delta(G)>(1-1/(r+4/3))|G|$, then $\chi(G)\le r+2$.
			\item If $\delta(G)>(1-1/(r+5/4))|G|$, then either $\chi(G)\le r+2$ or $G$ contains $K_{r-1}\vee \overline{C_7}$.
		\end{enumerate}
		For $r=1$, if $\delta(G)>5|G|/9$, then $G$ admits a homomorphism to $\overline{C_7}$. There are arbitrarily large locally bipartite graphs $G\nrightarrow \overline{C_7}$ whose minimum-degree ratios tend to $5/9$.
	\end{theorem}
	
	Combining the homomorphism conclusion of Theorem~\ref{thm:local} with a degree-counting argument yields the following corollary.
	\begin{corollary}\label{cor:rigid}
		If $G$ is locally bipartite, $\chi(G)>3$, and $\delta(G)\ge4|G|/7$, then $G\cong \overline{C_7}(t)$ for some positive integer $t$.
	\end{corollary}
	\begin{proof}
		Since $\delta(G)\ge4|G|/7>5|G|/9$, Theorem~\ref{thm:local} yields a homomorphism $\varphi:G\to\overline{C_7}$. Label the vertices of the original $C_7$ as $v_i$ for $i\in\mathbb{Z}_7$, in cyclic order, and put $V_i=\varphi^{-1}(v_i)$.
		Note that each $V_i$ is independent. Moreover, $V_i\ne\varnothing$ for every $i\in\mathbb{Z}_7$. Indeed, deleting any vertex from $\overline{C_7}$ leaves a 3-colourable graph: a perfect matching of the resulting path $P_6$ partitions its vertices into three independent pairs in the complement. If some $V_i$ were empty, $\varphi$ would therefore induce a 3-colouring of $G$, contrary to $\chi(G)>3$.
		
		For $x\in V_i$, the homomorphism property implies
		$4|G|/7\le d_G(x)\le |G|-|V_{i-1}|-|V_i|-|V_{i+1}|$.
		Thus $|V_{i-1}|+|V_i|+|V_{i+1}|\le3|G|/7$ for each $i\in\mathbb{Z}_7$. Summing over all seven indices gives
		\[|G|=\sum_{i\in \mathbb{Z}_7}|V_i|
		{{}={}} \frac13 \sum_{i\in \mathbb{Z}_7}(|V_{i-1}|+|V_i|+|V_{i+1}|)\le |G|.
		\] 
		Thus every triple has size $3|G|/7$. Subtracting the equations for consecutive indices gives $|V_{i-1}|=|V_{i+2}|$. Since $\gcd(3,7)=1$, it follows that for every $i\in \mathbb{Z}_7$,
		$|V_i|=|G|/7$.
		Let $t=|G|/7$.
		For each $x\in V_i$, all its neighbours lie in the four classes $V_j$ with $j-i\notin\{0,1,-1\}$, whose total size is $4t$. As $d_G(x)\ge4|G|/7=4t$, the vertex $x$ is adjacent to every vertex in these four classes. Hence $G\cong\overline{C_7}(t)$.
	\end{proof}
	
	\begin{lemma}\label{lem:common}
		Let $G$ be an $n$-vertex graph with minimum degree $\delta$. For an integer $r\ge1$ and an $r$-clique $R$ in $G$, let $M=\bigcap_{v\in V(R)}N_G(v)$. 
		If $\delta>(1-1/(r+c))n$, where $c>0$, then
		$|M|>\frac{c}{r+c}n$ and $\frac{\delta(G[M])}{|M|}>1-\frac1c$.
		The strict inequalities in the minimum-degree hypothesis and in both conclusions may be replaced simultaneously by non-strict inequalities.
	\end{lemma}
	
	\begin{proof}
		By the union bound, $|M|\ge\sum_{v\in V(R)}d_G(v)-(r-1)n$. Moreover, each vertex of $M$ has at most $n-|M|$ neighbours outside $M$. Therefore, whenever $M\ne\varnothing$,
		$|M|\ge r\delta-(r-1)n$ and $\delta(G[M])\ge\delta-(n-|M|)$.
		By the hypothesis $\delta>(1-1/(r+c))n$, we obtain
		\[
		|M|\ge r\delta-(r-1)n
		>r\left(1-\frac1{r+c}\right)n-(r-1)n
		=\frac{c}{r+c}n
		\]
		and
		\[
		\frac{\delta(G[M])}{|M|}
		\ge\frac{\delta-(n-|M|)}{|M|}
		=1-\frac{n-\delta}{|M|}
		>1-\frac{n/(r+c)}{cn/(r+c)}
		=1-\frac1c.
		\]
		If $\delta\ge(1-1/(r+c))n$, the same estimates give $|M|\ge cn/(r+c)$ and $\frac{\delta(G[M])}{|M|}
		\ge1-\frac1c$.
	\end{proof}

	For $r,t\in\mathbb{N}$, define
	\[
	G_r(t)=K_r(5t)\vee C_7(t)
	\qquad\text{and}\qquad
	G_r^{*}(t)=K_{r-1}(3t)\vee\overline{C_7}(t).
	\]
	When $r=1$, the complete multipartite graph in the second construction is empty, so $G_1^{*}(t)=\overline{C_7}(t)$.
	The graphs $G_r(t)$ and $G_r^{*}(t)$ are regular, with orders $(5r+7)t$ and $(3r+4)t$ and degrees $(5r+2)t$ and $(3r+1)t$, respectively. 
	Consequently,
	\[
	\frac{\delta(G_r(t))}{|G_r(t)|}=\frac{5r+2}{5r+7}
	\qquad\text{and}\qquad
	\frac{\delta(G_r^{*}(t))}{|G_r^{*}(t)|}=\frac{3r+1}{3r+4}.
	\]
	Moreover, we have
	\[
	\chi(G_r(t))=r+\chi(C_7)=r+3
	\qquad\text{and}\qquad
	\chi(G_r^{*}(t))=r-1+\chi(\overline{C_7})=r+3.
	\]
	Their forbidden-subgraph properties will be verified in the proof below.
	
	\begin{proof}[Proof of Theorem~\ref{thm:main}]
		Let $r\ge1$, $k\ge2$, and $n\ge1000(r+1)k^8$. Let $G$ be an $n$-vertex $H_{r,k}$-free graph with minimum degree $\delta=\delta(G)\ge\delta_{r,k}n$. Assume that $\chi(G)>r+2$. We first prove that $G$ has one of the two stated structures, and then verify that the constructions defined above establish sharpness.
		
		\medskip\noindent\textbf{Case 1: $k=2$.}
		
		Since $\delta\ge\delta_{r,2}n>\left(1-\frac1{r+4/3}\right)n$,
		Theorem~\ref{thm:local}(a) implies that $G$ is not $r$-locally bipartite. Hence there exists an $r$-clique $R$ such that $Q=G[N_G(R)]$ is non-bipartite. As $G$ is $H_{r,2}$-free, $Q$ is $C_5$-free.
		The non-strict form of Lemma~\ref{lem:common}, applied with $c=7/5$, gives $\delta(Q)\ge2|Q|/7$ and $|Q|\ge \frac{7n}{5r+7}
		\ge\frac{n}{r+1}
		\ge1000\cdot2^8$.
		Theorem~\ref{thm:cycle}(ii) therefore gives $Q\cong C_7(t)$ for some $t\in\mathbb{N}$. In particular, $|Q|=7t$ and $\delta(Q)=2t$. Thus,
		\[
		2t=\delta(Q)\ge\delta-(n-|Q|)
		\ge |Q|-\frac{n}{r+7/5}
		\ge |Q|-\frac{5|Q|}{7}=2t,
		\]
		where the bound $|Q|\ge \frac{7n}{5r+7}$ gives $n\le(5r+7)|Q|/7$. Equality holds throughout, so
		$ n=(5r+7)t$ and $\delta=(5r+2)t$.
		
		Let $U=V(G)\setminus V(Q)$, so $|U|=5rt$. Every vertex of $Q$ has exactly $2t$ neighbours in $Q$ and degree at least $(5r+2)t$ in $G$. It is therefore adjacent to all $5rt$ vertices of $U$.
		We claim that $G[U]$ is $K_{r+1}$-free. Otherwise, take $r$ vertices of such a clique and let $z$ be its remaining vertex. Four consecutive vertices of a seven-cycle in $Q$, together with $z$, form a $C_5$ in the common neighbourhood of the selected $r$ vertices. This gives a copy of $H_{r,2}$ in $G$, a contradiction.
		Furthermore,
		\begin{equation}
			\delta(G[U])\ge\delta-|Q|=5(r-1)t
			=\left(1-\frac1r\right)|U|.
		\end{equation}
		Since $r$ divides $|U|$, the minimum-degree bound and Tur\'an's theorem imply
		\[
		\left(1-\frac1r\right)\frac{|U|^2}{2}
		\le e(G[U])
		\le\left(1-\frac1r\right)\frac{|U|^2}{2}.
		\]
		The equality characterization in Tur\'an's theorem gives $G[U]\cong K_r(5t)$. This also covers $r=1$, when $G[U]$ is edgeless. Since all edges between $U$ and $V(Q)$ are present, $G\cong K_r(5t)\vee C_7(t)$.
		
		\medskip\noindent\textbf{Case 2: $k\ge3$.}
		
		The minimum-degree hypothesis is equivalent to
		$n-\delta\le\frac{3n}{3r+4}$.
		For any $r$-clique $R$, the non-strict form of Lemma~\ref{lem:common}, with $c=4/3$, shows that its common neighbourhood induces a graph of order at least
		$\frac{4n}{3r+4}\ge\frac{n}{r+1}\ge1000k^8$
		and minimum degree at least one quarter of its order. This graph is $C_{2k+1}$-free, since $G$ is $H_{r,k}$-free. As $\delta_k\le2/9<1/4$, Theorem~\ref{thm:cycle}(i) implies that it is bipartite. Thus $G$ is $r$-locally bipartite; this also holds if $G$ has no $r$-clique.
		
		If $r=1$, Corollary~\ref{cor:rigid} applies directly to $G$ and gives $G\cong\overline{C_7}(t)$. We therefore assume $r\ge2$. Since
		$\delta\ge\left(1-\frac1{r+4/3}\right)n
		>\left(1-\frac1{r+5/4}\right)n$
		and $\chi(G)>r+2$, Theorem~\ref{thm:local}(b) gives a copy of $K_{r-1}\vee\overline{C_7}$ in $G$. Let $R$ be the $(r-1)$-clique in this copy, and put $Q=G[N_G(R)]$. Then $Q$ contains $\overline{C_7}$ and hence $\chi(Q)\ge4$. For every $v\in V(Q)$, the graph induced by $N_Q(v)$ is the common-neighbourhood graph of the $r$-clique $R\cup\{v\}$ in $G$, and is bipartite. Therefore $Q$ is locally bipartite.
		
		By the common-neighbourhood estimate and the minimum-degree hypothesis, $|Q|\ge(r-1)\delta-(r-2)n\ge\frac{7n}{3r+4}$.
		Consequently,
		\begin{equation}
			\delta(Q)\ge\delta-(n-|Q|)
			\ge|Q|-\frac{3n}{3r+4}
			\ge\frac{4|Q|}{7}.\label{eq:long-core-degree}
		\end{equation}
		Since $Q$ is locally bipartite and $\chi(Q)\ge4$, Corollary~\ref{cor:rigid} implies that $Q\cong\overline{C_7}(t)$ for some $t\in\mathbb{N}$. Hence $|Q|=7t$ and $\delta(Q)=4t$. Substituting these values into \eqref{eq:long-core-degree} gives
		\[
		4t\ge\delta-n+7t
		\ge7t-\frac{3n}{3r+4}
		\ge4t.
		\]
		Equality therefore holds throughout. The last equality yields $n=(3r+4)t$, and the first gives $\delta=n-3t$. Thus
		$\delta=(3r+1)t$.
		
		Let $U=V(G)\setminus V(Q)$, so $|U|=3(r-1)t$. Every vertex of $Q$ has exactly $4t$ neighbours in $Q$ and at least $\delta=(3r+1)t$ neighbours in $G$. It is therefore adjacent to all vertices of $U$.
		If $G[U]$ contained an $r$-clique, then its common neighbourhood in $G$ would contain the non-bipartite graph $Q$, contrary to the $r$-local bipartiteness of $G$. Hence $G[U]$ is $K_r$-free. Moreover,
		\[
		\delta(G[U])\ge\delta-|Q|=3(r-2)t
		=\left(1-\frac1{r-1}\right)|U|.
		\]
		The minimum-degree bound and Tur\'an's theorem imply
		$e(G[U])=\left(1-\frac1{r-1}\right)\frac{|U|^2}{2}$.
		Since $|U|=3(r-1)t$, the equality characterization in Tur\'an's theorem gives $G[U]\cong K_{r-1}(3t)$. Together with the complete adjacency between $U$ and $V(Q)$, this yields $G\cong K_{r-1}(3t)\vee\overline{C_7}(t)$, proving~(ii).
		
		In both cases, a non-$(r+2)$-colourable graph satisfying the hypotheses has minimum degree exactly $\delta_{r,k}n$. In particular, $\delta(G)>\delta_{r,k}n$ guarantees $\chi(G)\le r+2$.
		It remains to verify the forbidden-subgraph properties of the two constructions defined above. In $G_r(t)$, an $r$-clique $R$ uses $j\in\{0,1,2\}$ vertices from $C_7(t)$. If $j=0$, its common-neighbourhood graph is $C_7(t)$, whose odd girth is seven. If $j=1$, the common-neighbourhood graph is the join of one unused outer class with the union of two cycle classes. This union is independent, so the common-neighbourhood graph is bipartite. If $j=2$, the common neighbourhood within $C_7(t)$ is empty, leaving two unused outer classes, which also induce a bipartite graph. Hence no $r$-clique has a $C_5$ in its common neighbourhood, and $G_r(t)$ is $H_{r,2}$-free.
		
		Every vertex neighbourhood in $\overline{C_7}$ is a $P_4$, so $\overline{C_7}$ is locally bipartite and $K_4$-free. An $r$-clique $R$ in $G_r^{*}(t)$ therefore uses $j\in\{1,2,3\}$ vertices from $\overline{C_7}(t)$. For $j=1$, its common-neighbourhood graph is a blow-up of $P_4$. For $j=2$, it is the join of one unused outer class with an independent common neighbourhood in $\overline{C_7}(t)$. For $j=3$, the common neighbourhood within $\overline{C_7}(t)$ is empty, leaving two unused outer classes. Each possible common-neighbourhood graph is bipartite. Thus $G_r^{*}(t)$ is $r$-locally bipartite and is $H_{r,k}$-free for every $k\ge1$.
		
		The degree and chromatic-number calculations above show that $G_r(t)$ attains $\delta_{r,2}$ and has chromatic number $r+3$, while $G_r^{*}(t)$ attains $\delta_{r,k}$ for every $k\ge3$ and also has chromatic number $r+3$. These constructions have arbitrarily large order. Together with the strict minimum-degree conclusion, they give $\chias(H_{r,k},r+2)=\delta_{r,k}$, completing the proof.
	\end{proof}
	
	\section{A fixed-target homomorphism threshold for wheels}\label{sec:fixed-target}
	
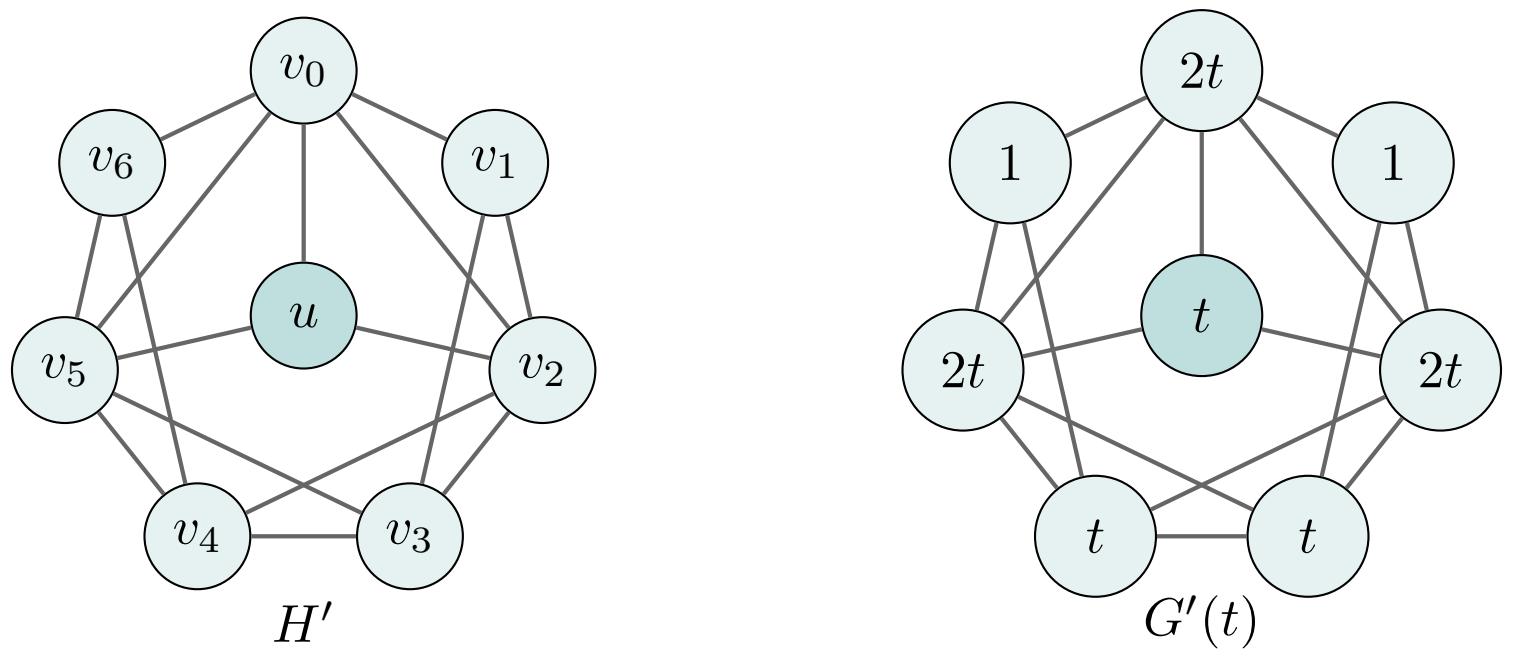
\begin{figure}[htbp]
		\centering
		\begin{tikzpicture}[scale=.9,every node/.style={font=\small}]
			\begin{scope}[xshift=-3.3cm]
				\foreach \i in {0,...,6}{
					\node[circle,draw,fill=teal!10,minimum size=7mm,inner sep=1pt]
					(a\i) at ({90-360*\i/7}:1.8) {$v_{\i}$};
				}
				\foreach \i/\j in {0/1,1/2,2/3,3/4,4/5,5/6,6/0,0/2,1/3,2/4,3/5,4/6,5/0}
				\draw[thick,gray!80!black] (a\i)--(a\j);
				\node[circle,draw,fill=teal!25,minimum size=7mm,inner sep=1pt]
				(au) at (0,0) {$u$};
				\foreach \i in {0,2,5}\draw[thick,gray!80!black] (au)--(a\i);
				\node at (0,-2.25) {$H'$};
			\end{scope}
			\begin{scope}[xshift=3.3cm]
				\foreach \i/\w in {0/2t,1/1,2/2t,3/t,4/t,5/2t,6/1}{
					\node[circle,draw,fill=teal!10,minimum size=8mm,inner sep=1pt]
					(b\i) at ({90-360*\i/7}:1.8) {$\w$};
				}
				\foreach \i/\j in {0/1,1/2,2/3,3/4,4/5,5/6,6/0,0/2,1/3,2/4,3/5,4/6,5/0}
				\draw[thick,gray!80!black] (b\i)--(b\j);
				\node[circle,draw,fill=teal!25,minimum size=8mm,inner sep=1pt]
				(bu) at (0,0) {$t$};
				\foreach \i in {0,2,5}\draw[thick,gray!80!black] (bu)--(b\i);
				\node at (0,-2.25) {$G'(t)$};
			\end{scope}
		\end{tikzpicture}
		\caption{The graph $H'$ (left) and its blow-up $G'(t)$ (right). In the right-hand diagram, each circle represents an independent set of the indicated size, and each line represents all edges between the corresponding sets.}
		\label{fig:five-ninths-counterexample}
	\end{figure}

	We use Illingworth's graph $H_2^+$ and the weighting in \cite[Fig.~2]{IlB}. In the following labelling, $H'=H_2^+$ has vertex set $\{v_0,\ldots,v_6,u\}$ and edge set
	\[
	E(H')=
	\bigl(\{v_iv_{i+1},v_iv_{i+2}:i\in\mathbb{Z}_7\}
	\setminus\{v_1v_6\}\bigr)
	\cup\{uv_0,uv_2,uv_5\},
	\]
	where subscripts are taken modulo $7$. For $t\in\mathbb{N}$, let $G'(t)$ be the blow-up of $H'$ whose independent classes corresponding to $(v_0,v_1,v_2,v_3,v_4,v_5,v_6,u)$ have sizes $(2t,1,2t,t,t,2t,1,t)$, respectively. The classes corresponding to $v_1$ and $v_6$ consist of vertices of degree $5t$, and every other vertex has degree at least $5t$. Consequently,
	$|G'(t)|=9t+2$, $\delta(G'(t))=5t$, and 
	$\lim_{t\to\infty}\frac{\delta(G'(t))}{|G'(t)|}=\frac59$.
	Figure~\ref{fig:five-ninths-counterexample} illustrates the base graph and the resulting construction.
	
	The graph $H'$ is locally bipartite and edge-maximal with this property \cite[Section~1.1]{IlB}. Its eight vertices have distinct neighbourhoods, so it is twin-free. By \cite[Lemma~4]{IlB}, a homomorphism $H'\to\overline{C_7}$ would imply that $H'$ is an induced subgraph of $\overline{C_7}$, which is impossible since $|H'|=8>7$. Blow-ups preserve local bipartiteness, and every class of $G'(t)$ is nonempty. Thus $G'(t)$ is locally bipartite and contains $H'$, so it is $W_{2k+1}$-free for every $k\ge1$ and admits no homomorphism to $\overline{C_7}$.

	\begin{proof}[Proof of Theorem~\ref{thm:fixed}]
		Set $d=\max\{\frac 59,\frac{2k+3}{4k+4}\}$.
		For $k\ge2$ and $n\ge2000k^8$, let $G$ be an $n$-vertex $W_{2k+1}$-free graph with $\delta(G)> dn$.
		For $x\in V(G)$, put $Q=G[N_G(x)]$. This graph is $C_{2k+1}$-free, since such a cycle together with $x$ would form a copy of $W_{2k+1}$.
		Since $|Q|\ge\delta(G)>dn$,
		\begin{equation}\nonumber
			\frac{\delta(Q)}{|Q|}
			\ge \frac{\delta(G)-(n-|Q|)}{|Q|}
			=1-\frac{n-\delta(G)}{|Q|}
			>  2-\frac1d.
		\end{equation}
		By the definition of $\delta_k$, $1/(2-\delta_k)=\max\{6/11,(2k+3)/(4k+4)\}\le d$. Thus $\delta(Q)>\delta_k|Q|$. Also $|Q|>dn>n/2\ge1000k^8$, so Theorem~\ref{thm:cycle}(i) implies that $Q$ is bipartite.
		Hence $G$ is locally bipartite. Since $\delta(G)>5n/9$, Theorem~\ref{thm:local} gives a homomorphism from $G$ to $\overline{C_7}$,
		proving the upper bound.
		
		For $k=2,3$, use the graphs $\overline{K}_{(2k+1)t}\vee C_{2k+3}(t)$, which have order $(4k+4)t$ and minimum degree $(2k+3)t$. These graphs are $W_{2k+1}$-free: the neighbourhood of a vertex in the independent part is $C_{2k+3}(t)$, which has odd girth $2k+3$, whereas the neighbourhood of a vertex in the cycle part is bipartite.
		Each graph contains $W_{2k+3}$. Under any homomorphism to $\overline{C_7}$, its rim would map to the bipartite neighbourhood of the image of its hub, which is impossible. Thus they do not map to $\overline{C_7}$, and their ratios give the lower bound.
		For $k\ge4$, the graphs $G'(t)$ defined above are $W_{2k+1}$-free and do not admit a homomorphism to $\overline{C_7}$. Their minimum-degree ratios tend to $5/9$, giving the required lower bound in this range.
	\end{proof}
	
	\begin{remark*}
		For each fixed $k\ge3$, every sufficiently large $W_{2k+1}$-free graph $G$ with $\delta_{\hom}(W_{2k+1},\overline{C_7})<\delta(G)/|G|\le4/7$ admits a homomorphism to $\overline{C_7}$. The upper endpoint is attained by $\overline{C_7}(t)$, which is not tripartite.
	\end{remark*}

	\section{Stability near the colourability thresholds}\label{sec:stability}
	
	For the case $k=2$, we use the following homomorphism theorem of H\"aggkvist and Jin~\cite{HJ}.
	
	\begin{theorem}[\cite{HJ}]\label{thm:odd-girth-seven}
		Let $G$ be an $n$-vertex graph of odd girth at least $7$. If $\delta(G)>n/4$, then $G$ admits a homomorphism to $C_7$.
	\end{theorem}
	
	To apply Theorem~\ref{thm:odd-girth-seven} to the $C_5$-free common neighbourhood of an $r$-clique, we first show that the minimum-degree condition also excludes triangles when its order is at least $120$.
	
	\begin{lemma}\label{lem:triangle-exclusion}
		Let $G$ be an $n$-vertex $C_5$-free graph with $n\ge120$ and $\delta(G)>n/4$. Then $G$ is triangle-free.
	\end{lemma}
	
	\begin{proof}
		We call an edge triangular if it belongs to a triangle. We first show that every vertex is incident with at most $10$ triangular edges. Suppose otherwise, and let $u$ be incident with at least $11$ such edges. Set
		$S=\{v\in N_G(u):uv\text{ belongs to a triangle}\}$.
		Then $|S|\ge11$ and $G[S]$ has no isolated vertices. Choose a maximum cut $(A,B)$ of $G[S]$ with $|A|\ge |B|$. Every vertex has a neighbour across the cut, since otherwise moving that vertex to the opposite side would increase the size of the cut. Since $|A|\ge\lceil |S|/2\rceil\ge6$, we may choose a set $X$ of six vertices from $A$, and for each $x\in X$ choose a neighbour across the cut. Let $Y$ be the set of chosen neighbours. Thus $X$ and $Y$ are disjoint subsets of $N_G(u)$, $|X|=6$, and $|Y|\le6$.
		Then $|X\cup Y\cup \{u\}|\le13$ and
		\[
		\sum_{x\in X}|N_G(x)\setminus (X\cup Y\cup \{u\})|\ge6(\delta(G)-|X\cup Y\cup \{u\}|+1)
		>n-| X\cup Y\cup \{u\}|.
		\]
		Consequently, two distinct vertices $x_1,x_2\in X$ have a common neighbour $z\notin X\cup Y\cup \{u\}$. If $y\in Y$ is a neighbour of $x_1$, then $x_1yux_2zx_1$ is a five-cycle, a contradiction. This proves the asserted bound on the number of triangular edges incident with a vertex.
		
		Suppose now that $G$ contains a triangle $T=x_1x_2x_3x_1$. Delete every triangular edge except the three edges of $T$, and denote the resulting spanning subgraph by $G'$. The graph $G'$ is $C_5$-free and has $T$ as its unique triangle. Every vertex loses at most $10$ incident edges. Moreover, each vertex of $T$ retains its two incident edges of $T$ and therefore loses at most $8$ edges. Thus $\delta(G')\ge\delta(G)-10$ and $d_{G'}(x_i)\ge\delta(G)-8$ for $i\in\{1,2,3\}$. Let $N_i=N_{G'}(x_i)\setminus V(T)$.
		These sets are pairwise disjoint, since a vertex in $N_i\cap N_j$ with $i\ne j$ would form a triangle with $x_i$ and $x_j$ distinct from $T$. Furthermore,
		$|N_i|=d_{G'}(x_i)-2\ge\delta(G)-10>0$.
		Choose $y_i\in N_i$ for each $i$, and let
		$M_i=N_{G'}(y_i)\setminus\{x_i\}$.
		Then $|M_i|=d_{G'}(y_i)-1\ge\delta(G)-11$.
		
		We claim that the seven sets $V(T),N_1,N_2,N_3,M_1,M_2,M_3$ are pairwise disjoint. The absence of triangles other than $T$ implies that $M_i\cap V(T)=\varnothing$ and $M_i\cap N_i=\varnothing$. If $i\ne j$ and $z\in M_i\cap N_j$, then $zy_ix_ix_{6-i-j}x_jz$ is a five-cycle, a contradiction. If $z\in M_i\cap M_j$, then $zy_ix_ix_jy_jz$ is a five-cycle. Both possibilities contradict the fact that $G'$ is $C_5$-free, proving the claim. It follows that
		\[
		n\ge |V(T)|+\sum_{i=1}^{3}|N_i|+\sum_{i=1}^{3}|M_i|\ge3+3(\delta(G)-10)+3(\delta(G)-11)>\frac{3n}{2}-60\ge n{.}
		\]
		This contradiction shows that $G$ is triangle-free.
	\end{proof}
	
	We now establish the stability result for $H_{r,2}$-free graphs using Theorem~\ref{thm:odd-girth-seven} and Lemma~\ref{lem:triangle-exclusion}.
	
	\begin{theorem}\label{thm:short-stability}
		Let $r\ge1$, $n\ge200(r+1)$, and $0\le\varepsilon<\frac1{(5r+7)(3r+4)}$.
		Suppose that $G$ is an $n$-vertex $H_{r,2}$-free graph with
		$\chi(G)>r+2$ and
		$\delta(G)\ge\left(\frac{5r+2}{5r+7}-\varepsilon\right)n$.
		There is a partition of $V(G)$ into nonempty parts $A_1,\ldots,A_r,V_0,\ldots,V_6$ such that
		\[
		\left(\frac5{5r+7}-\left(r+\frac25\right)\varepsilon\right)n
		\le |A_i|\le
		\left(\frac5{5r+7}+\varepsilon\right)n
		\qquad(1\le i\le r),
		\]
		and
		\[
		\left(\frac1{5r+7}-(3r+4)\varepsilon\right)n
		\le |V_i|\le
		\left(\frac1{5r+7}+(2r+3)\varepsilon\right)n
		\qquad(0\le i\le6).
		\]
		Let $J$ be the blow-up of $K_r\vee C_7$ on these parts, with $V_0,\ldots,V_6$ corresponding to the cycle vertices in cyclic order. There is a set $Z\subseteq\bigcup_{i=1}^r A_i$ with $|Z|\le(5r+7)\varepsilon n$ such that $G-Z\subseteq J-Z$ and
		$|E(G){{}\triangle{}} E(J)|\le2(5r+7)\varepsilon n^2$.
	\end{theorem}
	
	\begin{proof}
		Write $\delta=\delta(G)$. The bound on $\varepsilon$ gives $\delta>(1-3/(3r+4))n$. By Theorem~\ref{thm:local}(a), there is an $r$-clique $R$ with $V(R)=\{u_1,\ldots,u_r\}$ such that $Q=G[N_G(R)]$ is non-bipartite. Lemma~\ref{lem:common}, with $c=4/3$, gives
		\begin{equation}\label{eq:short-relative-degree}
			|Q|>\frac{4n}{3r+4}\ge200
			\qquad \text{and}\qquad\delta(Q)>\frac{|Q|}{4}.
		\end{equation}
		We also have
		\begin{align}
			|Q|\ge n-r(n-\delta(G))\qquad \text{and}\qquad
			n-\delta\le\left(\frac5{5r+7}+\varepsilon\right)n.\label{eq:short-degree-deficit}
		\end{align}
		Since $G$ is $H_{r,2}$-free, $Q$ is $C_5$-free. Lemma~\ref{lem:triangle-exclusion} therefore implies that $Q$ is triangle-free, so Theorem~\ref{thm:odd-girth-seven} gives a homomorphism $\varphi:Q\to C_7$.
		
		Let $C$ be a shortest odd cycle in $Q$. Then $|C|\ge7$ and $C$ is chordless. If a vertex outside $C$ has two neighbours $x,y$ on $C$, the odd $x$--$y$ path on $C$, together with these two incident edges, forms an odd cycle. By minimality of $C$, the complementary even path has length two. Thus any two such neighbours have distance two on $C$. Since $|C|\ge7$, there cannot be three such neighbours. Every vertex of $Q$ therefore has at most two neighbours on $C$, giving
		\[
		|C|\cdot\delta(Q)\le\sum_{x\in V(C)}d_Q(x)
		=\sum_{v\in V(Q)}|N_Q(v)\cap V(C)|\le2|Q|.
		\]
		By \eqref{eq:short-relative-degree}, $|C|<8$, so $|C|=7$. The restriction of $\varphi$ to $C$ is an isomorphism onto $C_7$: every proper subgraph of $C_7$ is bipartite, so the image of this odd cycle is all of $C_7$. Put $V_i=\varphi^{-1}(i)$ for $i\in\mathbb{Z}_7$, and choose four consecutive vertices $v_0,\ldots,v_3$ of $C$ with $v_i\in V_i$. All seven classes are nonempty.
		For $i\in\mathbb{Z}_7$ and $x\in V_i$, we have 
		\begin{equation}\label{eq:short-neighbour-sum}
			|V_{i-1}|+|V_{i+1}|
			\ge d_Q(x)
			\ge \delta(Q)
			\ge |Q|-n+\delta.
		\end{equation}
		Summing \eqref{eq:short-neighbour-sum} and applying \eqref{eq:short-degree-deficit} yields
		\begin{equation}\label{iq4.3}
			\begin{aligned}
				\sum_{i\in\mathbb{Z}_7}
				\bigl(|V_{i-1}|+|V_{i+1}|-|Q|+n-\delta\bigr)
				&=7(n-\delta)-5|Q|
				\le(5r+7)\varepsilon n.
			\end{aligned}
		\end{equation}
		Fix $i\in\mathbb{Z}_7$. Summing \eqref{eq:short-neighbour-sum} over $j\in\{i,i-3,i+3\}$ and over $j\in\{i-2,i-1,i+1,i+2\}$, respectively, gives
		$|Q|-|V_i|\ge3(|Q|-n+\delta)$ and  $|Q|+|V_i|\ge4(|Q|-n+\delta)$.
		Hence
		\begin{equation}\label{eq:short-fibre-sizes}
			3|Q|-4(n-\delta)\le |V_i|
			\le3(n-\delta)-2|Q|.
		\end{equation}
		Substituting \eqref{eq:short-degree-deficit} into \eqref{eq:short-fibre-sizes} gives the asserted bounds on $|V_i|$.
		
		Let $U=V(G)\setminus V(Q)$ and $Z_1=\{x\in U:|V(R)\setminus N_G(x)|\ge2\}$. Every vertex of $U$ misses at least one vertex of $R$. Counting nonadjacent pairs gives
		\[
		\begin{aligned}
			|Z_1|&\le\sum_{x\in U}\bigl(|V(R)\setminus N_G(x)|-1\bigr)
			=\sum_{j=1}^r(n-d_G(u_j))-|U|\le r(n-\delta)-n+|Q|.
		\end{aligned}
		\]
		Set $Z_2=\bigcup_{i=0}^3(U\setminus N_G(v_i))$. By \eqref{eq:short-neighbour-sum}, each of the seven expressions $|V_{i-1}|+|V_{i+1}|-|Q|+n-\delta$ is nonnegative. We may therefore extend the sum over $i=0,\ldots,3$ to all seven indices and apply \eqref{iq4.3} to obtain
		\[
		\begin{aligned}
			|Z_2|&\le\sum_{i=0}^3|U\setminus N_G(v_i)|
			=\sum_{i=0}^3\bigl(n-|Q|-d_G(v_i)+d_Q(v_i)\bigr)
			\le\sum_{i=0}^3\bigl(|V_{i-1}|+|V_{i+1}|-|Q|+n-\delta\bigr)\\
			&\le\sum_{i\in\mathbb{Z}_7}\bigl(|V_{i-1}|+|V_{i+1}|-|Q|+n-\delta\bigr)=7(n-\delta)-5|Q|.
		\end{aligned}
		\]
		Let $Z=Z_1\cup Z_2$. Combining the bounds on $|Z_1|$ and $|Z_2|$ with \eqref{eq:short-degree-deficit} gives
		\begin{equation}\label{eq:short-Z-bound}
			\begin{aligned}
				|Z|&\le |Z_1|+|Z_2|\le r(n-\delta)-n+|Q|+7(n-\delta)-5|Q|\\
				&=(r+7)(n-\delta)-n-4|Q|\le(5r+7)(n-\delta)-5n
				\le(5r+7)\varepsilon n.
			\end{aligned}
		\end{equation}
		
		Partition $U\setminus Z$ into $A'_1,\ldots,A'_r$ according to the unique non-neighbour $u_j$ in $R$. Each $A'_j$ contains $u_j$, and every vertex of $A'_j$ is adjacent to all vertices of $R\setminus\{u_j\}$ and $v_0,\ldots,v_3$. An edge $xy$ in $A'_j$ would yield the $r$-clique $(R\setminus\{u_j\})\cup\{x\}$ with the five-cycle $yv_0v_1v_2v_3y$ in its common neighbourhood. Thus $A'_j$ is independent and $|A'_j|\le n-\delta$.
		By \eqref{eq:short-degree-deficit},
		$\sum_{j=1}^r(n-\delta-|A'_j|)
		=r(n-\delta)-n+|Q|+|Z|\ge |Z|$.
		We can therefore distribute $Z$ among the $A'_j$ to obtain nonempty parts $A_j\supseteq A'_j$ of $U$ with $|A_j|\le n-\delta$. By \eqref{iq4.3},
		\[
		\begin{aligned}
			|A_j|&=n-|Q|-\sum_{i\ne j}|A_i|\ge n-|Q|-(r-1)(n-\delta)\ge n-\left(r+2/5\right)(n-\delta).
		\end{aligned}
		\]
		Together with \eqref{eq:short-degree-deficit} and $|A_j|\le n-\delta(G)$, this proves the asserted outer-part bounds.
		
		Let $J$ be the blow-up of $K_r\vee C_7$ on $A_1,\ldots,A_r,V_0,\ldots,V_6$. The independence of the $A'_j$ and the homomorphism $\varphi$ give $G-Z\subseteq J-Z$.
		Every vertex of $G-Z$ has degree at least $\delta-|Z|$. For $x\in V_i$, subtracting this from $d_{J-Z}(x)=n-|Z|-|Q|+|V_{i-1}|+|V_{i+1}|$ and applying \eqref{eq:short-neighbour-sum} and \eqref{iq4.3} gives
		\begin{equation}\label{eq:short-core-missing}
			\begin{aligned}
				d_{J-Z}(x)-d_{G-Z}(x)
				&\le |V_{i-1}|+|V_{i+1}|-|Q|+n-\delta\le(5r+7)\varepsilon n.
			\end{aligned}
		\end{equation}
		For $x\in A'_j$, the same subtraction gives at most $n-|A'_j|-\delta$ missing neighbours. Using \eqref{eq:short-degree-deficit}, \eqref{iq4.3} and \eqref{eq:short-Z-bound}, we obtain
		\begin{equation}\label{eq:short-outer-missing}
			\begin{aligned}
				n-\delta-|A'_j|
				&\le\sum_{i=1}^r(n-\delta-|A'_i|)=r(n-\delta)-n+|Q|+|Z|\\
				&\le\frac{(5r+7)(n-\delta(G))-5n}{5}+|Z|
				\le2(5r+7)\varepsilon n.
			\end{aligned}
		\end{equation}
		By \eqref{eq:short-core-missing} and \eqref{eq:short-outer-missing}, each vertex of $G-Z$ has at most $2(5r+7)\varepsilon n$ missing neighbours in $J-Z$. Counting each missing edge twice gives
		$e(J-Z)-e(G-Z)\le(5r+7)\varepsilon n^2$.
		Hence,
		\[
		|E(G){{}\triangle{}} E(J)|
		\le e(J-Z)-e(G-Z)+|Z|n
		\le2(5r+7)\varepsilon n^2.\qedhere
		\]
	\end{proof}
	
	For $k\ge3$, we first prove a stability theorem for $r$-locally bipartite graphs.
	
	\begin{theorem}\label{thm:local-stability}
		Let $r\ge1$ and
		$0\le\varepsilon<\frac1{(3r+4)(4r+5)}$.
		Suppose $G$ is an $n$-vertex $r$-locally bipartite graph with
		$\chi(G)>r+2$ and
		$ \delta(G)\ge\left(\frac{3r+1}{3r+4}-\varepsilon\right)n$.
		There is a partition into nonempty parts
		$V(G)=A_1\cup\cdots\cup A_{r-1}\cup V_0\cup\cdots\cup V_6$,
		where $V_0,\ldots,V_6$ are indexed modulo seven and the outer parts are absent when $r=1$. The part sizes satisfy
		\[
		\left(\frac3{3r+4}-\left(r+\frac13\right)\varepsilon\right)n
		\le |A_i|\le\left(\frac3{3r+4}+\varepsilon\right)n
		\qquad(1\le i\le r-1)
		\]
		and
		\[
		\left(\frac1{3r+4}-(r+1)\varepsilon\right)n
		\le |V_i|\le
		\left(\frac1{3r+4}+(2r+3)\varepsilon\right)n
		\qquad(0\le i\le 6).
		\]
		Let $J$ be the blow-up of $K_{r-1}\vee \overline{C_7}$ on this partition, with consecutive indices nonadjacent in $\overline{C_7}$.
		There is a set $Z\subseteq\bigcup_i A_i$, with $|Z|\le(3r+4)\varepsilon n$, such that $G-Z\subseteq J-Z$, and
		$|E(G){{}\triangle{}} E(J)|
		\le(3r+6)\varepsilon\left(1+\frac{3r+4}{3}\varepsilon\right)n^2$.
		In particular, when $r=1$, $Z=\varnothing$, $G\subseteq J$, and the sharper bound
		$ e(J)-e(G)\le\frac{\varepsilon+35\varepsilon^2}{2}\,n^2$
		holds. 
	\end{theorem}
	
	For the proof with $r\ge2$, we use a weighted version of the degree-deficit count.
	
	\begin{lemma}\label{lem:weighted}
		Let $r\ge2$, and let $G$ be an $r$-locally bipartite graph of order $n$ containing $R\vee \overline{C_7}$, where $R=\{u_1,\ldots,u_{r-1}\}$ is a clique. Write $Q=G[N_G(R)]$.
		There is a set $Z\subseteq V(G)\setminus V(Q)$ with
		$|Z|\le(3r+1)n-(3r+4)\delta(G)$
		such that
		$V(G)\setminus(V(Q)\cup Z)$ has a partition into independent sets $A'_1,\ldots,A'_{r-1}$. For each $i$, the set $A'_i$ contains $u_i$, and each of its vertices is adjacent to every vertex of $R\setminus\{u_i\}$ and to all seven vertices of the specified copy of $\overline{C_7}$.
	\end{lemma}
	\begin{proof}
		For each vertex $v$, define $\omega(v)=3r+1-3|N_G(v)\cap R|-|N_G(v)\cap V(\overline{C_7})|$. We verify that $\omega(v)\ge0$. If $v\notin V(Q)$, then $v$ is nonadjacent to at least one vertex of $R$, so $\omega(v)\ge 3r+1-3(r-2)-7=0$. If $v\in V(Q)$, then $Q$ is locally bipartite, as $G$ is $r$-locally bipartite. Since $\alpha(\overline{C_7})=2$, any five vertices of $\overline{C_7}$ induce a non-bipartite graph. Thus $v$ has at most four neighbours in the specified copy of $\overline{C_7}$, giving $\omega(v)\ge 3r+1-3(r-1)-4=0$. Summing over all vertices yields 
		\begin{equation}
			\begin{split}
				\sum_{v\in V(G)}\omega(v)
				&=(3r+1)n-3\sum_{i=1}^{r-1}d_G(u_i)
				-\sum_{b\in V(\overline{C_7})}d_G(b)\\
				&{\le(3r+1)n-(3r+4)\delta(G).}
			\end{split}\label{eq:weighted-deficit}
		\end{equation}
		In particular, the right-hand side is nonnegative. If $\omega(v)=0$, then $3|N_G(v)\cap R|+|N_G(v)\cap V(\overline{C_7})|=3r+1$. Hence either $v\in V(Q)$, or $v$ is nonadjacent to exactly one $u_i$ and adjacent to all seven vertices of $\overline{C_7}$. Assign each such vertex outside $Q$ to $A'_i$, where $u_i$ is its unique non-neighbour in $R$. Let $Z$ be the set of vertices outside $Q$ with positive weight. Each such vertex contributes at least one to the nonnegative sum in \eqref{eq:weighted-deficit}, giving the asserted bound on $|Z|$.
		
		By construction, $u_i\in A'_i$. If $x,y\in A'_i$ were adjacent, then $(R\setminus\{u_i\})\cup\{x,y\}$ would be an $r$-clique with $\overline{C_7}$ in its common neighbourhood, contrary to $r$-local bipartiteness. Thus each $A'_i$ is independent.
	\end{proof}
	
	\begin{proof}[Proof of Theorem~\ref{thm:local-stability}]
		Write $\delta=\delta(G)$.
		The assumed range of $\varepsilon$ gives
		\begin{equation}
			n-\delta
			\le\left(\frac3{3r+4}+\varepsilon\right)n
			<\frac{4n}{4r+5}\label{eq:nonneighbour-bound}
		\end{equation}
		and
		\[
		\frac{\delta}{n}>\frac{3r+1}{3r+4}-\frac{1}{(3r+4)(4r+5)}=\frac{4r+1}{4r+5}.
		\]
		When $r=1$, take $Q=G$. Then $Q$ is locally bipartite, $\chi(Q)>3$, and $\delta(Q)>5|Q|/9$.
		When $r\ge2$, Theorem~\ref{thm:local}(b) gives a copy of $R\vee \overline{C_7}$, where $R$ is an $(r-1)$-clique. Put $Q=G[N_G(R)]$. 
		We have
		\begin{equation}
			|Q|\ge n-(r-1)(n-\delta)\qquad\text{and}\qquad
			\delta(Q)\ge |Q|-(n-\delta).\label{eq:core-neighbourhood}
		\end{equation}
		Applying the strict form of Lemma~\ref{lem:common} to the $(r-1)$-clique $R$ with $c=9/4$ gives $\delta(Q)>5|Q|/9$. For each $v\in V(Q)$, the neighbourhood of $v$ in $Q$ is the common neighbourhood in $G$ of the $r$-clique $R\cup\{v\}$. Hence $Q$ is locally bipartite. Moreover, $Q$ contains $\overline{C_7}$, so $\chi(Q)\ge4$.
		
		In both cases, Theorem~\ref{thm:local} gives a homomorphism $Q\to\overline{C_7}$. Let $V_0,\ldots,V_6$ be its fibres, indexed modulo seven. All seven fibres are nonempty: otherwise $Q$ would map to a proper induced subgraph of $\overline{C_7}$, which is $3$-colourable, contrary to $\chi(Q)\ge4$.
		
		A vertex in $V_i$ has no neighbours in $V_{i-1}\cup V_i\cup V_{i+1}$. Then, for $0\le i\le6$,
		\begin{equation}\label{eq:triple-bound}
			|V_{i-1}\cup V_i\cup V_{i+1}|\le n-\delta.
		\end{equation}
		Each core vertex belongs to three of the sets in \eqref{eq:triple-bound}. Summing the differences between the right- and left-hand sides and using \eqref{eq:nonneighbour-bound} and \eqref{eq:core-neighbourhood}, we obtain
		\begin{equation}
			\begin{split}
				&\sum_{i\in\mathbb Z_7}\bigl(n-\delta-|V_{i-1}\cup V_i\cup V_{i+1}|\bigr)\\
				&=7(n-\delta)-3|Q|\le(3r+4)(n-\delta)-3n
				\le(3r+4)\varepsilon n.
			\end{split}\label{eq:degree-budget}
		\end{equation}
		Summing \eqref{eq:triple-bound} over the two indices $i-2$ and $i+2$ gives $|Q|-|V_i|\le2(n-\delta)$.
		The other five inequalities count each vertex of $V_i$ three times and every other core vertex twice. Consequently,
		$2|Q|+|V_i|\le5(n-\delta)$.
		We have proved
		\begin{equation}
			|Q|-2(n-\delta)\le |V_i|
			\le5(n-\delta)-2|Q|.\label{eq:fibre-bounds}
		\end{equation}
		Substituting the lower bound on $|Q|$ from \eqref{eq:core-neighbourhood} into both bounds gives
		\[
		|V_i|\ge n-(r+1)(n-\delta)
		\ge\left(\frac1{3r+4}-(r+1)\varepsilon\right)n,
		\]
		and
		\[
		|V_i|\le(2r+3)(n-\delta)-2n
		\le\left(\frac1{3r+4}+(2r+3)\varepsilon\right)n.
		\]
		This proves the asserted bounds on $|V_i|$.
		
		For $r=1$, there are no outer vertices. Take $Z=\varnothing$; the homomorphism gives $G\subseteq J$ for the blow-up $J$ of $\overline{C_7}$ on $V_0,\ldots,V_6$.
		Now suppose $r\ge2$.
		Apply Lemma~\ref{lem:weighted} to obtain independent sets $A'_1,\ldots,A'_{r-1}$ and an exceptional set $Z$ outside $Q$, with
		$|Z|\le(3r+1)n-(3r+4)\delta\le(3r+4)\varepsilon n$. Since each $A'_j$ is independent and nonempty, $|A'_j|\le n-\delta$. Moreover,
		\begin{equation}
			\begin{split}
				\sum_{j=1}^{r-1}(n-\delta-|A'_j|)
				&=(r-1)(n-\delta)-n+|Q|+|Z|\ge |Z|,
			\end{split}\label{eq:outer-capacity}
		\end{equation}
		where the last inequality is \eqref{eq:core-neighbourhood}.
		We may therefore distribute the vertices of $Z$ among the sets $A'_j$ so that the enlarged parts $A_j$ satisfy $|A_j|\le n-\delta$.
		Their total size is $n-|Q|$. Bounding each of the other $r-2$ outer parts by $n-\delta$ gives
		\[
		\begin{aligned}
			|A_j|
			&\ge n-|Q|-(r-2)(n-\delta)\ge n-\left(r+\frac13\right)(n-\delta)\\
			&{\ge\left(\frac3{3r+4}-\left(r+\frac13\right)\varepsilon\right)n,}
		\end{aligned}
		\]
		where the second inequality uses $3|Q|\le7(n-\delta)$, obtained by summing \eqref{eq:triple-bound}.
		The upper bound follows from \eqref{eq:nonneighbour-bound}. This proves  the asserted outer-part bounds.
		
		Let $J$ be the blow-up of $K_{r-1}\vee\overline{C_7}$ on $A_1\cup \cdots\cup A_{r-1}\cup V_0\cup \cdots\cup V_{6}$. After deleting $Z$, its outer parts are $A'_j=A_j\setminus Z$, each independent in $G$. The homomorphism $Q\to\overline{C_7}$ accounts for all edges within $Q$. Hence $G-Z\subseteq J-Z$.
		The following count applies to every $r\ge1$, with the outer sum empty for $r=1$. 
		
		Let
		$ M=e(J-Z)-e(G-Z)$.
		We count, at each vertex, the neighbours allowed by $J-Z$ that are missing in $G-Z$. Every vertex has degree at least $\delta-|Z|$ in $G-Z$.
		For a vertex in $V_i$, the number of allowed neighbours is
		$n-|Z|-(|V_{i-1}|+|V_i|+|V_{i+1}|)$.
		It therefore has at most
		$n-\delta-(|V_{i-1}|+|V_i|+|V_{i+1}|)$
		missing neighbours.
		For a vertex in $A'_j$, the allowed neighbours are the vertices of $G-Z$ outside $A'_j$. It therefore has at most $n-\delta-|A'_j|$ missing neighbours. Summing over the core and outer parts gives
		\begin{equation}
			\begin{split}
				2M&\le
				\sum_{i=0}^6 |V_i|
				\bigl(n-\delta-(|V_{i-1}|+|V_i|+|V_{i+1}|)\bigr)+\sum_{j=1}^{r-1}|A'_j|\bigl(n-\delta-|A'_j|\bigr)\\
				&\le\bigl(5(n-\delta)-2|Q|\bigr)\bigl(7(n-\delta)-3|Q|\bigr)+(n-\delta)\bigl((r-1)(n-\delta)-n+|Q|+|Z|\bigr)\\
				&=\frac{n-\delta}{3}\bigl[(3r+1)n-(3r+4)\delta\bigr]+(n-\delta)|Z|+\frac23\bigl(7(n-\delta)-3|Q|\bigr)^2,
			\end{split}\label{eq:missing-count}
		\end{equation}
		where the second inequality uses \eqref{eq:fibre-bounds}, \eqref{eq:outer-capacity} and $|A'_j|\le n-\delta$; when $r=1$, the outer term vanishes. By Lemma~\ref{lem:weighted}, $|Z|\le(3r+4)\varepsilon n$ for $r\ge2$, and this bound also holds for $r=1$ since $Z=\varnothing$. Substituting this bound, \eqref{eq:nonneighbour-bound} and \eqref{eq:degree-budget} into \eqref{eq:missing-count} yields
		\[
		\begin{aligned}
			2M
			&\le\frac13\left(\frac3{3r+4}+\varepsilon\right)(3r+4)\varepsilon n^2+\left(\frac3{3r+4}+\varepsilon\right)(3r+4)\varepsilon n^2
				+\frac23\bigl((3r+4)\varepsilon n\bigr)^2\\
			&=4\varepsilon n^2+2(r+2)(3r+4)\varepsilon^2 n^2.
		\end{aligned}
		\]
		It follows that
		$M\le2\varepsilon n^2+(r+2)(3r+4)\varepsilon^2 n^2$.
		
		For $r\ge 2$, all remaining discrepancies between $G$ and $J$ concern pairs with at least one endpoint in $Z$. There are at most $|Z|n$ such pairs. Consequently,
		\[
		\begin{aligned}
			|E(G)\triangle E(J)|
			&\le M+|Z|n\le(3r+6)\varepsilon n^2+(r+2)(3r+4)\varepsilon^2 n^2
			= (3r+6)\varepsilon\left(1+\frac{3r+4}{3}\varepsilon\right)n^2.
		\end{aligned}
		\]
		When $r=1$, $Q=G$, $Z=\varnothing$, and $M=e(J)-e(G)$. The first sum in \eqref{eq:missing-count} alone gives the stronger estimate
		\[
		\begin{aligned}
			2M
			&\le\max_i|V_i|\,\bigl(7(n-\delta)-3n\bigr)\le\left(\frac17+5\varepsilon\right)n\cdot7\varepsilon n
			=(\varepsilon+35\varepsilon^2)n^2.
		\end{aligned}
		\]
		This completes the proof.
	\end{proof}
	
	\begin{proof}[Proof of Theorem~\ref{thm:stability-summary}]
		{The case $k=2$ follows from Theorem~\ref{thm:short-stability}. Now let $k\ge3$, and suppose that $G$ satisfies the hypotheses of Theorem~\ref{thm:stability-summary}.}
		It suffices to show that $G$ is $r$-locally bipartite.
		Write $\delta=\delta(G)$. 
		The bound on $\varepsilon$ gives $\delta>(7r+2)n/(7r+9)$ when $k=3$ and $\delta>(4r+1)n/(4r+5)$ when $k\ge4$.
		For each $r$-clique $R$, let $Q=G[N_G(R)]$. 
		Applying Lemma~\ref{lem:common} with $c=9/7$ when $k=3$ and $c=5/4$ when $k\ge4$, we obtain $\delta(Q)/|Q|>2/9=\delta_3$ in the first case and $\delta(Q)/|Q|>1/5>\delta_k$ in the second.
		In both cases $n-\delta<n/(r+1)$, so $|Q|\ge n-r(n-\delta)>n/(r+1)\ge1000k^8$. Since $G$ is $H_{r,k}$-free, $Q$ is $C_{2k+1}$-free. Theorem~\ref{thm:cycle}(i) therefore implies that $Q$ is bipartite.
		The assumed bound on $\varepsilon$ also implies $\varepsilon<1/[(3r+4)(4r+5)]$. Applying Theorem~\ref{thm:local-stability} therefore gives all the conclusions of Theorem~\ref{thm:stability-summary} for $k\ge3$.
	\end{proof}

	%
	%
	%
	%

\end{document}